\documentclass[11pt]{amsart}
\usepackage{amsmath, amsthm}

\theoremstyle{definition}

\theoremstyle{remark}

\theoremstyle{definition}

\newtheorem{thm}{Theorem}[section]
\newtheorem{lem}[thm]{Lemma}

\newtheorem*{thmA}{Theorem~A}

\newtheorem{defn}[thm]{Definition}
\newtheorem{rem}[thm]{Remark}

\newtheorem{notation}[thm]{Notation}

\newcommand{\field}[1]{\mathbb{#1}}

\newcommand{\reals}{\ensuremath{\field{R}}}

\numberwithin{equation}{section}

\begin{document}

\title{Convergence of mean curvature flows with surgery}

\author{Joseph Lauer}
\address{Department of Mathematics, Yale University, New Haven, CT 06510}
%\curraddr{Department of Mathematics and Statistics,
%Case Western Reserve University, Cleveland, Ohio 43403}
\email{joseph.lauer@yale.edu}

%    General info

\date{September 5, 2011}


\begin{abstract}
Huisken and Sinestrari ~\cite{HS} have recently defined a surgery
process for mean curvature flow when the initial data is a
two-convex hypersurface in $\reals^{n+1}$ ($n\geq 3$). The process
depends on a parameter $H$. Its role is to initiate a surgery when
the maximum of the mean curvature of the evolving hypersurface
becomes $H$, and to control the scale at which each surgery is
performed.  We prove that as $H\to\infty$ the surgery process
converges to level set flow ~\cite{CCG91}~\cite{ES91}.
\end{abstract}

\maketitle

\section*{Introduction}

Huisken and Sinestrari ~\cite{HS} have recently defined a mean
curvature flow with surgery when the initial data is a two-convex
hypersurface in $\reals^{n+1}$ when $n\geq 3$.  The process depends
on a parameter $H$ ($H_3$ in the notation of ~\cite{HS}), which
controls both the maximal mean curvature and the scale at which each
surgery is performed.  In this note we investigate to what extent
the process depends on this parameter.

Recall that a smooth one-parameter family of hypersurface immersions
$F_t:M\to\reals^{n+1}$ is a solution to mean curvature flow if

$$
\frac{\partial F}{\partial t}(x,t)=\vec{H}(F(x,t)),
$$
where $\vec{H}$ is the mean curvature vector.  The first results
were obtained by Huisken ~\cite{Hu84} who proved that if the initial
data is convex and $n\geq 2$, then the mean curvature flow shrinks
the hypersurface to a round point.  The analogous result for curves
in the plane $(n=1)$ was proved by Gage and Hamilton~\cite{GH86},
and shortly after Grayson ~\cite{G87} showed that any embedded curve
in the plane evolves to become convex.  This means that the
classification of singularities is particularly simple for embedded
plane curves. However, when $n>1$ Grayson's Theorem no longer holds
and singularities other than round points may occur. The existence
of such a singularity was first proved rigourously by
Grayson~\cite{G89}, who gave the example of a barbell-like surface
which develops a neck-pinch.

As an evolving hypersurface becomes singular the maximum of the mean
curvature is unbounded, and hence constructing a surgery procedure
requires detailed information about the geometry of the hypersurface
in regions of high curvature. In the two-convex case, Huisken and
Sinestrari prove that such regions are diffeomorphic to $S^n$ or
$S^{n-1}\times S^1$, and are discarded during surgery, or are
neck-like regions in which the surgery replaces a topological
cylinder by a pair of convex disks. As the parameter $H$ increases
the surgeries are performed closer to the singular time and on
quantitatively thinner necks. The detailed estimates in ~\cite{HS}
controlling the length and width of the necks allow us to prove:

\begin{thmA}
As $H\to\infty$ the Huisken-Sinestrari surgery converges to level
set flow.
\end{thmA}

Since the limit is unique this result can be interpreted as a
stability theorem for level set flow.  Our approach is to use a
barrier argument:  We prove that for any $\epsilon>0$ there exists
$H>0$ so that the mean curvature flow with surgery performed with
parameter $H$ is disjoint (in space-time) from the level set flow of
the initial hypersurface shifted backwards in time by $\epsilon$.

Since the Ricci flow with surgery constructed for 3-manifolds (see~\cite{KL08} and ~\cite{P03}) 
also depends on a parameter, it is possible to consider the same question there.
One obstacle in this direction is that there is no natural candidate for the limiting object. 

{\bf Acknowledgements:} The author wishes to thank his advisor Bruce Kleiner for his guidance during work on the project,
and for suggesting the problem

%%%%%%%%%%%%%%%%%%%%%%%%%%%%%%%%%%%
%%%%% Weak notions of MCF %%%%%%%%%%%%%%%%%%%%
%%%%%%%%%%%%%%%%%%%%%%%%%%%%%%%%%%%%

\section{Weak notions of mean curvature flow}

In this section we recall (see ~\cite{HS}~\cite{I94}~\cite{W00}) two
ways in which the evolution of a smooth hypersurface can be extended
beyond a singularity: Level set flow and mean curvature flow with
surgery.

\begin{defn}  [Weak Set Flow] Let $K\subset\reals^{n+1}$ be
closed, and let $\{K_t\}_{t\geq 0}$ be a one-parameter family of
closed sets with initial condition $K_0=K$ such that the space-time
track $\cup (K_t\times\{t\})\subset\reals^{n+2}$ is closed.  Then
$\{K_t\}_{t\geq 0}$ is weak set flow for $K$ if for every smooth
mean curvature flow $\Sigma_t$ defined on $[a,b]$ we have
$$
K_a\cap \Sigma_a=\emptyset\Longrightarrow K_t\cap\Sigma_t=\emptyset
$$
for each $t\in[a,b]$.
\end{defn}

It is essentially the definition that weak set flows avoid smooth
mean curvature flows when the initial conditions are disjoint but a
stronger statement is true:  The distance between a weak set flow
and a smooth mean curvature flow is nondecreasing in $t$. Otherwise
one could translate the initial data in space and obtain a
contradiction to the definition of a weak set flow.

\begin{defn} [Level Set Flow]  The level set flow of a compact set
$K\subset\reals^{n+1}$, denoted $LSF(K)$, is the maximal weak set
flow. That is, a one-parameter family of closed sets $K_t$ with
$K_0=K$ such that if $\hat{K}_t$ is any weak set flow with
$\hat{K}_0=K$ then $\hat{K}_t\subset K_t$ for each $t\geq 0$.
\end{defn}

The existence of a maximal weak set flow is verified by taking the
closure of the union of all weak set flows with a given initial
data.  If $K_t$ is the weak set flow of $K$, we denote by~$\widehat{K}$ 
by the space-time track swept out by $K_t$.  That is, 
$$
\widehat{K}=\bigcup_{t\geq 0} K_t\times \{t\}\subset\reals^{n+2}.
$$

%We denote by LSF$(K)_T=K_T\times\{T\}$ the $t=T$ time-slice of
%LSF$(K)$, which is often identified with a subset or $\reals

The level set flow was introduced independently by Evans and Spruck
~\cite{ES91} and Chen, Giga and Goto ~\cite{CCG91}.   It was first
formulated in terms of viscosity solutions of partial differential
equations whereas the geometric definition above was first used by
Ilmanen ~\cite{I94}.

Another approach to constructing weak solutions to geometric
evolution equations has been to use a surgery procedure.  This idea
was first used by Hamilton~\cite{Ha83} to avoid the development of
singularities in Ricci Flow.

\begin{defn}
[Surgery, ~\cite{HS}] A mean curvature flow with surgery consists of
the following data: \newline 1) An initial smooth hypersurface
$\Sigma\subset\reals^{n+1}$.  \newline 2) Constants
$\omega_1<\omega_2<1$ and $H>0$.
\newline 3) A finite collection of
times $0<t_1<t_2\ldots <t_m$ called surgery times
 (let $t_0=0$).
 \newline 4) A collection of mean curvature flows $\Sigma_t^i$ on
$[t_i,t_{i+1}]$, with $\Sigma_0^0=\Sigma$, such that for each $i$
the maximum mean curvature on $\Sigma_t^i$ is $H$ and is achieved
only when $t=t_{i+1}$.
\newline 5) A surgery algorithm that consists of two steps:
\newline\indent i) At each surgery time a finite number of necks with mean
curvature greater than $\omega_1H$ are removed from
$\Sigma_{t_{i+1}}^{i}$ and replaced with convex caps with mean
curvature bounded by $\omega_2H$.  The operation of replacing a
single neck with two convex caps is called a standard surgery.
\newline\indent ii) Finitely many components of the hypersurface
constructed in i) are removed. These components are recognized as
being diffeomorphic to either $S^{n-1}\times S^1$ or
$S^n$.\newline\indent The result of the surgery algorithm is a
smooth hypersurface $\Sigma_{t_{i+1}}^{i+1}$ with mean curvature
bounded by $\omega_2H$.

We denote by $\Sigma_H\subset\reals^{n+2}$ the space-time track
swept-out by the hypersurfaces, and say that $\Sigma_H$ is a mean
curvature flow with surgery performed with parameter $H$.
\end{defn}

The main result of ~\cite{HS} is that a mean curvature flow with
surgery can be constructed when the initial data is a closed
two-convex hypersurface of dimension at least three.  A hypersurface
$\Sigma\subset\reals^{n+1}$ is two-convex if the sum of the two
smallest principal curvatures is everywhere nonnegative.  It is
proved that for any such initial data there exist $\omega_1$,
$\omega_2$ and $H_0>0$ so that the surgery may be performed with any
parameter $H\geq H_0$.  In particular, $\omega_1$ and $\omega_2$ can
be fixed independently of $H$.  It is also shown that if the initial
data is embedded then the hypersurface remains embedded even after a
surgery time.

It will be convenient to work with the regions bounded by the
evolving hypersurface.  Let $K\subset\reals^{n+1}$ be a compact
domain such that $\partial K$ is a smooth two-convex hypersurface.
Then if $\partial K_H$ is a mean curvature flow with surgery we
define $K_H\subset\reals^{n+2}$ to be the region of space-time such that
the $t=T$ time-slice of $K_H$ is the compact domain bounded 
by~$(\partial K_H)_T$.  The hypersurface $(\partial K_H)_t$ may not be
connected after the first surgery time.  However, the domains
bounded by the connected components of $(\partial K_H)_t$ will be
disjoint so that $(K_H)_t$ is well-defined.  Thus $K_H$ is an
evolution of a union of domains whose boundary is a mean curvature
flow with surgery performed with parameter $H$ in the sense defined
above.  We will also refer to $K_H$ as a mean curvature flow with
surgery.

If $K$ is a compact domain and $K_H$ is a mean curvature flow with
surgery constructed as in ~\cite{HS} then it is easy to verify that
$K_H$ is a weak set flow for $K$.   Note that this is not true if we
consider only the evolving hypersurfaces i.e., $\partial K_H$ is not
a weak set flow of $\partial K$.

\begin{notation} If $K_H$ is a mean curvature flow with surgery, and $T$ is a
surgery time, then we use $(\partial K_H)_T^-$ and $(\partial
K_H)_T^+$ to refer to the pre- and post-surgery hypersurfaces at
time $T$, and $(K_H)_T^-$ and $(K_H)_T^+$ to refer to the regions
they bound.
\end{notation}

%%%%%%%%%%%%%%%%%%%%%%%%%%%%%%%%
%%%% Convergence %%%%%%%%%%%%%%%%%%%%%
%%%%%%%%%%%%%%%%%%%%%%%%%%%%%%%%

\section{Convergence}

In this section we prove the convergence to level set flow.  Recall that $\widehat{K}$ denotes the space-time track of the level set flow of $K$.

\begin{thmA}\label{main} Let $K\subset\reals^{n+1}$, $n\geq 3$, be a compact
domain with $\partial K$ a smooth embedded two-convex hypersurface.  For $H$
sufficiently large let $K_H\subset\reals^{n+2}$ be the result of the
Huisken-Sinestrari surgery performed with parameter $H$, and initial
condition $(K_H)_0=K$.  Then

$$
\lim_{H\to\infty} K_H=\widehat{K}.
$$
\end{thmA}

\begin{rem} Convergence is with respect to the Hausdorff topology on
closed sets of $\reals^{n+2}$.
\end{rem}

Theorem $A$ follows from the following lemma regarding the surgery procedure, and a barrier argument.  As usual, $B_\epsilon(x)\subset\reals^{n+1}$ represents the ball of radius $\epsilon$ centered at $x$.

\begin{lem} \label{neck} Given $\epsilon>0$ there exists $H_0>0$ such that if
$H\geq H_0$, $T$ is a surgery time, and $x\in\reals^{n+1}$, then
$$
B_\epsilon(x)\subset (K_H)_T^-\Longrightarrow B_\epsilon(x)\subset
(K_H)_T^+.
$$
\end{lem}

The proof of Lemma ~\ref{neck} requires geometric information
regarding the necks along which a surgery is performed.  The
parameter $H$ here corresponds to $H_3$ in ~\cite{HS}, and
$\omega_1, \omega_2$ are the constants appearing in Definition 2.3.
Define $H_1=\omega_1 H$ and $H_2=\omega_2 H$. Furthermore,
$\epsilon_0, k, \Lambda$ are constant defined in ~\cite{HS} and depend only
on the initial hypersurface.

\proof [Proof of Lemma~\ref{neck}] Let $K_H$ be a mean curvature flow with surgery.

Since $T$ is a surgery time the Huisken-Sinestrari algorithm identifies 
a finite collection of subsets, $\{A_i\}_{i=1}^m$, which
cover the regions of $(\partial K_H)_T^-$ with mean curvature
greater than $H_2$.  There are three possibilities for the structure
of each $A_i$ depending on whether it has 0,1 or 2 boundary
components.

If $\partial A_i\neq\emptyset$ then for each component of $\partial
A_i$ a standard surgery is performed.  According to ~\cite{HS} there
exists an embedding $N:S^{n-1}\times[a,b]\to A_i$ with strong
geometric properties.  In particular, each
$\Sigma_z=N(S^{n-1}\times{z})$ has constant mean curvature
$\frac{n-1}{r_z}$, where $r_z$ is called the mean radius of
$\Sigma_z$.  If $\partial A_i$ consists of two
connected components then the map $N$ is a diffeomorphism.  In
general, $\partial A_i$ contains at least one of $\Sigma_a$ or
$\Sigma_b$ and the mean curvature on $\partial A_i$ is
$\frac{H_1}{2}$.

Suppose $\Sigma_a\subset\partial A_i$.  We consider the standard
surgery corresponding to $\Sigma_a$. Let $z_0\in[a,b]$ be the point
closest to $a$ such that the mean curvature on $\Sigma_{z_0}$ is
$H_1$.  The slice $\Sigma_{z_0}$ is sufficiently far from $\partial
A_i$ in the sense that $a<z_0-4\Lambda<z_0+4\Lambda<b$, where
$\Lambda\geq 10$. For simplicity we will assume that $z_0=0$. The
map $N$ can be extended (after first restricting it to
$S^{n-1}\times[-4\Lambda,4\Lambda]$) to a local diffeomorphism

$$
G: B_1^n\times[-4\Lambda,4\Lambda]\to\reals^{n+1}
$$
which is $\epsilon_0$-close in the $C^{k+1}$-norm to the standard
isometric embedding of some tube $B^n\times[-4\Lambda,4\Lambda]$ in
$\reals^{n+1}$ [~\cite{HS}, Prop. 3.25].  The standard surgery
removes $N(S^{n-1}\times[-3\Lambda,3\Lambda])$ and replaces it by
two convex caps contained in $G(B_1^n\times[-3\Lambda,3\Lambda])$,
and the result is again a smooth embedded hypersurface [~\cite{HS},
Thm. 3.26].  By the Jordan-Brouwer Separation Theorem for
hypersurfaces it follows that if $x\in (K_H)_T^-\setminus
G(B_1^n\times[-3\Lambda,3\Lambda])$ then $x$ will remain in the
interior of the hypersurface after the standard surgery.

Since $G$ is $\epsilon_0$-close to a standard tube and $\Lambda\geq
10$ is sufficiently large compared to $\epsilon_0$ we can choose
$H_0$ large enough (and hence the radius of the tube small enough)
so that if $H\geq H_0$ then
$$
B_\epsilon(x)\subset (K_H)_T^-\Longrightarrow B_\epsilon(x)\cap
G(B_1^n\times[-3\Lambda,3\Lambda])=\emptyset.
$$

With $H_0$ chosen in this way it follows that if $B_\epsilon(x)\subset
(K_H)_T^-$, then $B_\epsilon(x)$ lies in the region bounded by the hypersurface after a
standard surgery.  At each surgery time a finite number of
standard surgeries may be performed. However, the solid tubes
associated to the surgeries are disjoint and so the surgeries do not
interact.

It remains to verify that components discarded by 5)ii) of
Definition 2.3 do not bound a ball of radius $\epsilon$.  There are
three ways in which such a component can arise:
\newline 1) If $\partial A_i=\emptyset$ then $A_i$ is diffeomorphic to $S^n$ or
$S^{n-1}\times S^1$ and is discarded.
\newline 2) If $\partial A_i$ consists of a single
component then $A_i$ is homeomorphic to a ball. This corresponds to
the case where the curvature does not decrease significantly in one
direction along the neck. In this case only one standard surgery is
performed. After the standard surgery, the end of the cylinder with
high curvature will have become diffeomorphic to $S^n$ and will be
discarded.
\newline 3) If $\partial A_i$ consists of two components then a standard surgery is
performed for each boundary component and the result is two capped
cylinders and a component diffeomorphic to $S^2$.  The $S^2$
component is discarded.

In each case the construction in ~\cite{HS} guarantees that the mean
curvature of the component being removed is bounded from below by
$\frac{H_1}{2}$. Suppose $\Sigma$ is such a hypersurface, that $x$
lies in the region bounded by $\Sigma$ and that
$d=d(x,\Sigma)\geq\epsilon$.  If $y\in\Sigma$ realizes $d(x,\Sigma)$
then the mean curvature at $y$ is not more than $\frac{n}{d}\leq
\frac{n}{\epsilon}$ since $\Sigma$ $\cap$ int$(B_d(x))=\emptyset$.
This is a contradiction as long as
$H_0\geq\frac{2n}{\epsilon\omega_1}$.\qed

%%%%%%%%%%%%%%%%%%%%%%%%%%%%%%%%%%%%%%%%%%%%%%%%%%%
%%%%%%%%%%%%%%%%%%%%%%%%%%%%%%%%%%%%%%%%%%%%%%%%%%%

\proof [Proof of Theorem $A$] Given $\epsilon>0$ sufficiently small let $t_\epsilon>0$ be the time
such that
$$
d(\partial K,\partial K_{t_\epsilon})=\epsilon.
$$
Such a time exists since $\partial K$ is two-convex.  Let $\Omega_\epsilon\subset\reals^{n+2}$ 
be the level set flow $K_{t_e}$.  Then $\Omega_\epsilon$ is the level set flow of $K$ shifted 
backwards in time by $t_\epsilon$ (ignoring $t<0$). 

%and hence $\lim_{\epsilon\to 0}\Omega_\epsilon=\widehat{K}$.  

%Let  can be thought of as LSF$(K)$
%shifted backwards in time by $t_\epsilon$.

Let $H_0=H_0(\epsilon)$ be chosen as in Lemma~\ref{neck}.  

\smallskip

{\bf Claim:}  $\Omega_{\epsilon}\subset K_H$ for all $H\geq H_0$.

\smallskip

Let $T$ be the first surgery time of $K_H$.  Since $\partial K_H$ is
a smooth mean curvature flow on $[0,T)$ and $\Omega_\epsilon$ is a
weak set flow the distance between the two is nondecreasing on that
interval.  Thus $d((\Omega_\epsilon)_T,
(\partial K_H)_T^-)\geq\epsilon$ since $t_\epsilon$ was chosen so that $d((\Omega_\epsilon)_0,
(\partial K_H)_0)=\epsilon$.  
Applying Lemma ~\ref{neck} we conclude that $d((\Omega_\epsilon)_T, (\partial K_H)_T^+)\geq\epsilon$.  Since $(\partial
K_H)_T^+$ is a smooth hypersurface the argument can be repeated for
each of the subsequent surgery times. This proves the claim.

Since $\lim_{\epsilon\to 0}\Omega_\epsilon=\widehat{K}$ the claim implies that 
$\widehat{K}\subset$ lim$_{H\to\infty}K_H$ since
the limit of closed sets is closed.  Finally, since each mean
curvature flow with surgery is a weak set flow for $K$ the limit is
also and thus lim$_{H\to\infty}K_H\subset\widehat{K}$.\qed

%%%%%%%

%-------------------------------------------------------------------

\bibliographystyle{amsplain}

\end{document}